\input ./style/arxiv-general.cfg
\documentclass[MSNbibl,number,citesort,dvips]{arxbj}
\makeatletter
   \@ifpackageloaded{graphicx}{}{\usepackage{graphicx}}
\makeatother
\usepackage{dcolumn}

\volume{22}
\issue{2}
\pubyear{2016}
\firstpage{711}
\lastpage{731}
\doi{10.3150/14-BEJ673} % kopijuoti is laisko
\docsubty{FLA}

\makeatletter
\newcommand{\rrvert}{\vert}
\newcommand{\rrVert}{\Vert}
\newcommand{\llvert}{\vert}
\newcommand{\llVert}{\Vert}
\newtheorem{prop}{Proposition}[section]
\newtheorem{lem}{Lemma}[section]
\newtheorem{corol}{Corollary}[section]
\newremark{remark}{Remark}[section]
\newremark{example}{Example}[section]

\DeclareSymbolFont{AMSb}{U}{msb}{M}{n}
\newcommand{\Eout}{\mathbb{E}}
\newcommand{\PR}{\mathbb{P}}
\newcommand{\N}{\mathbb{N}}
\newcommand{\R}{\mathbb{R}}
\newcommand{\bX}{\mathbb{X}}
\newcommand{\bY}{\mathbb{Y}}
\newcommand{\bZ}{\mathbb{Z}}
\newcommand{\cd}{{\mathcal D}}
\newcommand{\cf}{{\mathcal F}}
\newcommand{\cl}{{\mathcal L}}
\newcommand{\cp}{{\mathcal P}}
\newcommand{\cq}{{\mathcal Q}}
\newcommand{\cv}{{\mathcal V}}
\newcommand{\cx}{{\mathcal X}}
\newcommand{\cdal}{{\cd_{\alpha}}}
\newcommand{\cdmo}{{\cd_{-1}}}
\newcommand{\cdpo}{{\cd_{+1}}}

\newcommand{\Falp}{F_\alpha}
\newcommand{\Fmalp}{F_{-\alpha}}
\newcommand{\Gamalp}{\Gamma_\alpha}
\newcommand{\Gamalptil}{\tilde{\Gamma}_\alpha}
\newcommand{\ximalp}{\xi_{-\alpha}}
\newcommand{\Upalp}{\Upsilon_\alpha}
\newcommand{\Pihat}{\hat{\Pi}}

\newcommand{\atil}{{\tilde{a}}}
\newcommand{\btil}{{\tilde{b}}}
\newcommand{\nablatil}{{\tilde{\nabla}}}
\newcommand{\Gtil}{{\tilde{G}}}
\newcommand{\Mtil}{\tilde{M}}
\newcommand{\Phitil}{\tilde{\Phi}}
\newcommand{\util}{{\tilde{u}}}
\newcommand{\vtil}{{\tilde{v}}}
\newcommand{\wtil}{{\tilde{w}}}
\newcommand{\cvtil}{{\tilde{\mathcal V}}}
\newcommand{\imathtil}{{\tilde{\imath}}}
\newcommand{\rhobar}{{\bar{\rho}}}

\newcommand{\bfa}{\mathbf{a}}
\newcommand{\bfatil}{\tilde{\mathbf{a}}}
\newcommand{\bfbtil}{\tilde{\mathbf{b}}}
\newcommand{\bfe}{\mathbf{e}}
\newcommand{\bfut}{\tilde{\mathbf{u}}}
\newcommand{\bfU}{\mathbf{U}}
\newcommand{\bfy}{\mathbf{y}}
\newcommand{\bfz}{\mathbf{z}}
\newcommand{\bfvt}{{\tilde{\mathbf{v}}}}
\newcommand{\bfV}{\mathbf{V}}
\newcommand{\bfW}{\mathbf{W}}

\newcommand{\E}{\mathbf{E}}
\newcommand{\Emu}{\mathbf{E}_\mu}
\newcommand{\EP}{\mathbf{E}_P}
\newcommand{\fndot}{\cdot}
\newcommand{\cond}{\mid}
\makeatother

\begin{document}
\begin{frontmatter}

\title{Infinite-dimensional statistical manifolds based on a balanced chart}
\runtitle{Infinite-dimensional manifolds}

\begin{aug}
%%%% inicialai - be tarpu
% Corresponding author: Nigel Newton - njn@essex.ac.uk% Updated by VTEXPTS2LaTeX.exe, 17.11.2014 09:13
\author[A]{\inits{N.J.}\fnms{Nigel J.}~\snm{Newton}\corref{}\thanksref{A}\ead[label=e1]{njn@essex.ac.uk}}
%%\runauthor{} %% auto
%\dedicated{}
\address[A]{School of Computer Science and Electronic Engineering,
University of Essex,
Wivenhoe Park, Colchester, CO4 3SQ, UK. \printead{e1}.}
\end{aug}

% HISTORY:
\received{\smonth{9} \syear{2013}}
\revised{\smonth{8} \syear{2014}}

% ABSTRACT
%
\begin{abstract}
We develop a family of infinite-dimensional Banach manifolds of
measures on an abstract
measurable space, employing charts that are ``balanced'' between the
density and
log-density functions. The manifolds, $(\tilde{M}_\lambda, \lambda\in
[2,\infty))$, retain
many of the features of finite-dimensional information geometry; in
particular, the
$\alpha$-divergences are of class $C^{\lceil\lambda\rceil-1}$,
enabling the definition
of the Fisher metric and $\alpha$-derivatives of particular classes of
vector fields.
Manifolds of \emph{probability} measures, $(M_\lambda, \lambda\in
[2,\infty))$, based on
centred versions of the charts are shown to be $C^{\lceil\lambda
\rceil-1}$-embedded
submanifolds of the $\tilde{M}_\lambda$. The Fisher metric is a
pseudo-Riemannian metric
on $\tilde{M}_\lambda$. However, when restricted to finite-dimensional embedded
submanifolds it becomes a Riemannian metric, allowing the full
development of the
geometry of $\alpha$-covariant derivatives. $\tilde{M}_\lambda$ and
$M_\lambda$ provide
natural settings for the study and comparison of approximations to posterior
distributions in problems of Bayesian estimation.
\end{abstract}

% KEYWORDS
% visi is mazosios raides ir pagal abecele
%
\end{frontmatter}

%s1 #&#
\section{Introduction} \label{seintro}

This paper develops a family of infinite-dimensional manifolds of
measures, each
containing a smoothly embedded submanifold of \emph{probability}
measures. It was
motivated by problems of Bayesian estimation, in which posterior
distributions have to
be computed from a variety of partial observations. This can rarely be
done exactly
owing to issues of dimension and nonlinearity, and the study of
approximations is
contingent on the development of appropriate measures of error. The
manifolds we
construct have metrics suited to such problems.

Suppose, for example, that $X\dvtx \Omega\rightarrow\bX$ and $Y\dvtx \Omega
\rightarrow\bY$
are random variables defined on a common probability space $(\Omega,\cf,\PR)$, taking
values in metric spaces $\bX$ and $\bY$, respectively. $X$~is the
\emph{estimand}
whose posterior distribution we seek, and $Y$ is the observable. Let
$\cp$ be
the set of probability measures on the Borel subsets of $\bX$. Under
mild conditions
(see, e.g., \cite{lish1}) an abstract Bayes formula defines a
regular conditional
distribution for $X$ given $Y$, $\Pi\dvtx \bY\rightarrow\cp$. (For any
Borel set
$B\subseteq X$, $\Pi(\fndot)(B)\dvtx \bY\rightarrow[0,1]$ is measurable, and
$\PR(X\in B\cond Y)=\Pi(Y)(B)$.) In the applications we have in mind,
$\Pi(\bY)$
is typically of infinite dimension and so we need to construct
approximations of
the form $\Pihat\dvtx \bY\rightarrow\cq\subset\cp$, where $\cq$ is of
finite dimension.

Single estimation objectives, such as minimum mean-square error in the
approximation of
a real-valued random variable $f(X)$, induce their own specific
measures of
error on $\cp$, but these may not be easy to use. On the other hand,
if $f$ is
sufficiently regular, then a more generic measure of error such as the
$L^2$ metric
on densities may be useful. If $\mu\in\cp$ is a reference measure
with respect
to which $\Pi(y)$ and $\Pihat(y)$ have densities $\pi(y)$ and $\hat
{\pi}(y)$, then the
difference between the minimum mean-square error estimate of $f(X)$ and the
mean of $f$ under $\Pihat(y)$ can be bounded by means of
the Cauchy--Schwarz
inequality:
%
%e1 #&#
%
\begin{equation}
(\E_{\Pi(y)}f-\E_{\Pihat(y)}f )^2 \le\Emu f^2
\Emu\bigl(\pi(y)-\hat{\pi}(y)\bigr)^2. \label{eql2bnd}
\end{equation}
Although, in this context, the $L^2$ metric on densities induces an appropriate
topology on $\cp$, it may still be poor in practice. This is so, for example,
if $f$ is the indicator function of a rare, but important, event.
Moreover, we often
need generic measures of error that are suitable for a \emph{variety} of
objectives.
This is especially important if the underlying estimation problem is inherently
multi-objective, as is the case, for example, when tracking the
movement of many
objects.

The mean-square error of $\E_{\Pihat(Y)} f$ admits the orthogonal
decomposition:
\[
\Eout \bigl(f(X)-\E_{\Pihat(Y)}f \bigr)^2 = \Eout
\E_{\Pi(Y)} (f-\E_{\Pi(Y)}f )^2 + \Eout (
\E_{\Pi(Y)}f-\E_{\Pihat(Y)}f )^2.
\]
The first term on the right-hand side is the \emph{estimation error}
arising from
the limitations of the observation $Y$; the second term is the \emph{approximation error}
arising from the use of $\Pihat$ instead of $\Pi$. When comparing
errors for more than
one random variable, it is natural to normalise the approximation errors
by their associated estimation errors---there is no point in
approximating the
conditional mean, $\E_{\Pi(y)}f$, with high precision if it is itself
a poor estimate
of $f(X)$. With this in mind, we might propose the following extreme,
multi-objective,
mean-square measure of error on $\cp$:
%
%e2 #&#
%
\begin{eqnarray}
\cd(Q\cond P) &:= & \sup_{f\in\cl^2(P)}\frac{(\EP f-\E_Qf)^2}{\EP(f-\EP f)^2}
\nonumber
\\
& = & \sup_{f\in F} \bigl(\EP f(1-\mathrm{d}Q/\mathrm{d}P) \bigr)^2
\label{eqchisq}
\\
& = & \EP(1-\mathrm{d}Q/\mathrm{d}P)^2,
\nonumber
\end{eqnarray}
where $\cl^2(P)=\{f\dvtx \bX\rightarrow\R\dvtx  \EP f^2<\infty\}$ and $F$ is
the subset of such
functions having zero mean and unit variance. This is the $\chi^2$-divergence.
Although extreme, it illustrates a feature of many multi-objective
measures of error:
they ensure that probabilities of events that are small are
approximated with greater
absolute accuracy than those that are large. The $L^p$ metrics on
densities fail in
this respect. (A related disadvantage is that spaces of probability
densities have
boundaries, which can create problems with numerical methods.) A
commonly used, less
extreme, multi-objective measure of error is the \emph{Kullback--Leibler}
divergence.
This is widely used in variational Bayesian estimation. (See, e.g.,
\cite{grim1,smqu1}.)

The regularity of the Kullback--Leibler divergence was central to the
design of the
manifolds in this paper. Each manifold is covered by a single chart,
which places its
elements in one-to-one correspondence with those of a Banach space.
Because of this,
the manifolds are also \emph{metric spaces} of measures, with metrics
tailored (at least
locally) to problems of Bayesian estimation. The manifolds are large
enough to include
exact posterior distributions in many problems. They also include, as
smoothly embedded
submanifolds, a large variety of finite-dimensional families of
probability measures,
on which approximations can be based.

The study and approximation of \emph{nonlinear filters} (an application pursued
elsewhere by the author) was a particular motivation. A nonlinear
filter computes the
posterior distributions of a Markov \emph{signal process} from randomly-perturbed
observations that become progressively available in time. For a modern
perspective on
the theory and application of nonlinear filtering, the reader is
referred to
\cite{crro1}. Approximations based on information geometric
projections onto
finite-dimensional exponential families were studied in \cite{bhlg1}.

The equations of nonlinear filtering are often expressed in terms of the
``un-normalised'' version of the posterior distribution obtained when
the marginal
density of the observation is omitted from the denominator in Bayes'
formula. This
satisfies the so-called \emph{Zakai equation}, which has a particularly
simple (bilinear)
form. A manifold of finite measures with a suitable metric is a natural
space for such
un-normalised posteriors. We develop a family of such manifolds in
Section~\ref{semanifmt}, not only because of this application, but
also because many of the
properties of the \emph{statistical} manifolds of Section~\ref
{semanifm} are best
understood in the context of their embedding in these larger manifolds.
The manifolds
are also natural settings in which to study and compare
finite-dimensional statistical
manifolds that admit the full geometry of $\alpha$-covariant derivatives.

The paper is structured as follows. Section~\ref{seinfgeo} provides a brief
introduction to information geometry. Section~\ref{semanifmt}
introduces the
one-parameter family of manifolds of finite measures,
$((\Mtil_\lambda,\tilde{\phi}_\lambda), \lambda\in[2,\infty
))$, and studies on them the
properties of Amari's $\alpha$-embedding maps. Section~\ref
{semanifm} develops the
family of manifolds of \emph{probability} measures
$((M_\lambda,\phi_\lambda), \lambda\in[2,\infty))$ in which the
chart $\phi_\lambda$ is
a ``centred'' version of $\tilde{\phi}_\lambda$. Section~\ref
{sestatdiv} studies the
properties of the $\alpha$-divergences on $\Mtil_\lambda$ and
$M_\lambda$, defining the
Fisher metric, and a limited notion of $\alpha$-parallel transport on
the tangent
bundle. Some examples of finite-dimensional embedded submanifolds of
$\Mtil_\lambda$
and $M_\lambda$ are outlined in Section~\ref{sefindim}. A sketch of
some of the
results of Sections~\ref{semanifm} and \ref{sealphm} was given,
without proofs, in
\cite{newt6}.

%s2 #&#
\section{Information geometry} \label{seinfgeo}

We begin by reviewing a classical finite-dimensional example: the exponential
statistical manifold. (See, e.g., \cite{amar1}.) Let $(\bX,\cx,\mu)$ be a
probability space supporting real-valued random variables $(\eta_i;
i=1,\ldots,d)$
with the following properties: (i) the variables
$(1,\eta_1,\eta_2,\ldots,\eta_d)$ are linearly independent elements
of $L^0(\mu)$,
that is, $\mu(\alpha+\sum_i y^i\eta_i = 0)=1$ if and only if
$\alpha=0$ and $\R^d\ni y=0$;
(ii) $\Emu\exp(\sum_iy^i\eta_i)<\infty$ for all $y$ in a non-empty
open subset
$B\subseteq\R^d$. For each $y\in B$, let $P_y$ be the probability
measure on $\cx$
with density
%
%e3 #&#
%
\begin{equation}
\frac{\mathrm{d}P_y}{\mathrm{d}\mu} = \exp \biggl(\sum_iy^i
\eta_i-c(y) \biggr), \label{eqexpman}
\end{equation}
where $c(y)=\log\Emu\exp(\sum_iy^i\eta_i)$, and let $N:=\{P_y\dvt
y\in B\}$. It follows
from (i) that the map $B\ni y\mapsto P_y\in N$ is a bijection. Let
$\theta\dvtx N\rightarrow B$ be its inverse; then $(N,B,\theta)$ is an
\emph{exponential
statistical manifold}, with an atlas comprising the single chart
$\theta$. We can
think of a \emph{tangent vector} at $P\in N$, $U$, as being an
equivalence class of
differentiable curves passing through $P$: two curves (expressed in
coordinates),
$(\bfy(t)\in B, t\in(-\varepsilon,\varepsilon))$ and
$(\bfz(t)\in B, t\in(-\varepsilon,\varepsilon))$, being equivalent at
$P$ if
$\bfy(0)=\bfz(0)=\theta(P)$ and $\dot{\bfy}(0)=\dot{\bfz}(0)$.
The \emph{tangent space}
at $P$, $T_PN$, is the linear space of all such tangent vectors, and is
spanned by the
vectors $(\partial_i; i=1,\ldots,d)$, where $\partial_i$ is the
equivalence class
containing the curve $(\bfy_i(t):=\theta(P)+t\bfe_i, t\in
(-\varepsilon,\varepsilon))$,
and $\bfe_i^j$ is equal to the Kr\"{o}necker delta. The \emph{tangent
bundle} is the
disjoint union $TN:= \bigcup_{P\in N}(P,T_PN)$, and admits the global chart
$\Theta\dvtx TN\rightarrow B\times\R^d$, where
$\Theta^{-1}(y,u)=(\theta^{-1}(y),u^i\partial_i)$. If a function
$f\dvtx N\rightarrow\R^n$
is differentiable, and $U\in T_PN$, then we write
%
%e4 #&#
%
\begin{equation}
Uf = u^i\partial_i f:= u^i\frac{\mathrm{d}}{\mathrm{d}t}
\bigl(f\circ\theta^{-1}\bigr) \bigl(\bfy_i(t)\bigr)\biggm|
_{t=0} = u^i\frac{\partial(f\circ\theta^{-1})}{\partial y^i}(y),
\end{equation}
where $(y,u)=\Theta(P,U)=(\theta(P),U\theta)$, and we have used the
Einstein summation
convention, that indices appearing once as a superscript and once as a
subscript are
summed out.

For each $\alpha\in[-1,1]$, let $\cdal\dvtx N\times N\rightarrow
[0,\infty)$ be the
$\alpha$-divergence
%
%e5 #&#
%
\begin{equation}
\cdal(P\mid Q):= \cases{ \displaystyle \E_Q\frac{\mathrm{d}P}{\mathrm{d}Q}\log
\frac{\mathrm{d}P}{\mathrm{d}Q}, &\quad if $\alpha=-1$,
\vspace*{3pt}\cr
\displaystyle\frac{4}{1-\alpha^2}
\biggl(1-\E_Q \biggl(\frac{\mathrm{d}P}{\mathrm{d}Q} \biggr)^{(1-\alpha)/2}
\biggr), &\quad if $\alpha\in(-1,1)$,
\vspace*{3pt}\cr
\displaystyle \E_Q\log
\frac{\mathrm{d}Q}{\mathrm{d}P}, &\quad if $\alpha=1$.} \label{eqalphdiv}
\end{equation}
(The Kullback--Leibler divergence corresponds to the case $\alpha=-1$.)
These are
of class $C^\infty$; their mixed second derivatives define the \emph{Fisher metric} as
a Riemannian metric on $N$: for any $P\in N$, any $U,V\in T_PN$, and
any $\alpha\in[-1,1]$,
%
%e6 #&#
%
\begin{equation}
\langle U, V \rangle_P:= -UV\cdal= u^ig(P)_{i,j}v^j,
\label{eqexpfish}
\end{equation}
where $U$ and $V$ act on the first and second argument of $\cdal$,
respectively, and
%
%e7 #&#
%
\begin{equation}
g(P)_{i,j}:= \langle\partial_i, \partial_j
\rangle_P = \EP(\eta_i-\EP\eta_i) (
\eta_j-\EP\eta_j). \label{eqfishmat}
\end{equation}
The mixed third derivatives of the $\alpha$-divergences define a
family of
\emph{covariant derivatives} on~$N$. If~$\bfU,\bfV\dvtx N\rightarrow TN$
are sufficiently
smooth vector fields of $N$ then the Chentsov--Amari \mbox{$\alpha$-}covariant
derivative is
defined as follows:
%
%e8 #&#
%
\begin{equation}
\nabla_\bfU^\alpha\bfV(P) = \bfU{\mathbf v}^k(P)
\partial_k + \Gamma_\alpha(P)_{i,j}^k{\mathbf
u}(P)^i{\mathbf v}(P)^j\partial_k. \label{eqcovder}
\end{equation}
Here ${\mathbf u}(P)=\bfU(P)\theta$, ${\mathbf v}(P)=\bfV(P)\theta$, and
the \emph{Christoffel symbols}
are as follows
%
%e9 #&#
%
\begin{eqnarray}\label{eqchrissym}
\Gamma_\alpha(P)_{i,j}^k & = & -
g(P)^{k,l}\partial_i\partial_j
\partial_l\cdal
\nonumber\\[-8pt]\\[-8pt]\nonumber
& = & \frac{1-\alpha}{2}g(P)^{k,l} \EP(\eta_i-\EP
\eta_i) (\eta_j-\EP\eta_j) (
\eta_l-\EP\eta_l),
\nonumber
\end{eqnarray}
where $g(P)^{k,l}$ is the $(k,l)$ element of the inverse of the matrix $g(P)$,
$\partial_i$ and $\partial_j$ act on the first argument of $\cdal$, and
$\partial_l$ acts on the second argument.

The covariant derivatives $\nabla^\alpha$ and $\nabla^{-\alpha}$
are \emph{dual} in the
sense that, for appropriately smooth vector fields $\bfU,\bfV$ and
$\bfW$,
%
%e10 #&#
%
\begin{equation}
\bfU\langle\bfV, \bfW\rangle_P = \bigl\langle\nabla_\bfU^\alpha
\bfV, \bfW \bigr\rangle_P + \bigl\langle\bfV, \nabla_\bfU^{-\alpha}
\bfW \bigr\rangle_P. \label{eqexpcodu}
\end{equation}
Each $\alpha$-covariant derivative defines a notion of \emph{parallel
transport} on the
tangent bundle. Equation (\ref{eqexpcodu}) shows that, if two tangent
vectors at
base point $P$ are parallel transported along a differentiable curve,
one according to
$\nabla^\alpha$ and the other according to $\nabla^{-\alpha}$, then
their Fisher dot
product remains constant. The $\alpha$-covariant derivatives thus
generalise the
Levi--Civita covariant derivative of Riemannian geometry, which
corresponds to the
special case $\alpha=0$.

Information geometry is the study of such structures, and has a history
going back to
the work of Rao \cite{rao1}. It derives its importance from the
fundamental role
played by the Fisher information in estimation theory. An example
application in
asymptotic statistics is to decompose the error of a \emph{second-order efficient
estimator} into a term arising from the choice of the estimator and
terms arising from
the curvature of the parametric model from which the estimate is
chosen. (See
Chapter~4 in \cite{amar1}.) For more applications, from a variety of
fields, the
reader is referred to~\cite{niba1}.

The literature on information geometry is dominated by the study of
finite-dimensional
manifolds of probability measures (parametric models) such as
$(N,B,\theta)$ above.
See \cite{amba1,amar1,barn1,chen1,muri1} and the references
therein for further
information. However, these are not always sufficiently inclusive for
the Bayesian
applications outlined in Section~\ref{seintro}, and any extension of
the ideas to the
non-parametric case must be based on charts with respect to which the
$\alpha$-divergences are suitably smooth. As is clear from the first
equation in
(\ref{eqalphdiv}), the smoothness properties of the Kullback--Leibler
divergence are
closely connected with those of the density, $\mathrm{d}P/\mathrm{d}Q$, and its log
(considered as
elements of dual function spaces). In the series of papers
\cite{cepi1,gipi1,piro1,pise1}, G.~Pistone and his coworkers developed an
infinite-dimensional exponential statistical manifold on an abstract
probability space
$(\bX,\cx,\mu)$. (See, also, \cite{gras1,zhha1}.) Probability
measures in the
manifold are mutually absolutely continuous with respect to the
reference measure
$\mu$, and the manifold is covered by the charts $s_Q(P)=\log \mathrm{d}P/\mathrm{d}Q-\E
_Q\log \mathrm{d}P/\mathrm{d}Q$ for
different ``patch-centric'' probability measures $Q$. These readily
give $\log \mathrm{d}P/\mathrm{d}Q$
the desired regularity, but require ranges that are subsets of
exponential Orlicz
spaces in order to do the same for $\mathrm{d}P/\mathrm{d}Q$. The exponential Orlicz
manifold is the
natural infinite-dimensional extension of the exponential manifold
$(N,B,\theta)$
described above; it has a strong topology, under which the $\alpha
$-divergences are of
class $C^\infty$. Variants of the Chentsov--Amari covariant derivatives
are defined on
it in \cite{gipi1}. However, with the exception of the case $\alpha
=1$, they are not
defined on the tangent bundle. If $\alpha\in(-1,1)$, for example, the
$\alpha$-connection is defined on the vector bundle whose fibre at
base point $P$ is
the Lebesgue space $L^{2/(1-\alpha)}(P)$.

This approach is highly inclusive, but is technically demanding and
leads to manifolds
that are larger than needed in many applications. The author's aim in
\cite{newt4} and
the present paper was to construct simpler infinite-dimensional
statistical manifolds
appropriate to problems in Bayesian estimation. The manifolds we
construct differ from
one another in the numbers of derivatives that the $\alpha
$-divergences admit.
A minimal requirement is a mixed second derivative since this is needed
in the
construction of the Fisher metric. It is achieved in a Hilbert setting
in \cite{newt4}.
However, it is also useful for the $\alpha$-divergences to admit
higher derivatives
so that notions of parallel transport can be developed. This is
achieved here in the
context of Banach manifolds.

%s3 #&#
\section{The manifolds of finite measures} \label{semanifmt}

Let $(\bX, \cx, \mu)$ be a probability space. For some $\lambda\in
[2,\infty)$, we
consider the set, $\Mtil$ ($=\Mtil_\lambda$), of finite measures on
$\cx$ satisfying
the following conditions:
\begin{enumerate}[(M3)]
\item[(M1)] $P$ is mutually absolutely continuous with respect to $\mu$;
\item[(M2)] $\Emu p^\lambda< \infty$;
\item[(M3)] $\Emu\llvert  \log p\rrvert  ^\lambda< \infty$.
\end{enumerate}
(We denote measures in $\Mtil$ by the upper-case letters $P$, $Q$,
\ldots, and their
densities with respect to $\mu$ by the corresponding lower case
letters, $p$, $q,\ldots.$) In order to control both the density $p$ and its log, we
employ a ``balanced''
chart involving their sum. Let $\Gtil$ ($=\Gtil_\lambda:=L^\lambda
(\mu)$) be the
Lebesgue space of (equivalence classes of) random variables $\atil\dvtx \bX
\rightarrow\R$
for which $\Emu\llvert  \atil\rrvert  ^\lambda<\infty$, and let $\tilde{\phi
}\dvtx \Mtil\rightarrow\Gtil$ be
defined as follows:
%
%e11 #&#
%
\begin{equation}
\tilde{\phi}(P) = p-1 + \log p. \label{eqphitdef}
\end{equation}

%pr3.1 #&#
%
\begin{prop} \label{prbiject}
$\tilde{\phi}$ is a bijection onto $\Gtil$.
\end{prop}

\begin{pf}
For $y\in(0,\infty)$ let $\theta(y) = y-1+\log y$; then $\inf_y
\theta(y)=-\infty$,
$\sup_y \theta(y)=+\infty$, and $\theta$ is of class $C^\infty$
with first derivative
$\theta^{(1)}(y) = 1 + y^{-1} > 0$. So, according to the inverse
function theorem,
$\theta\dvtx (0,\infty)\rightarrow\R$ is a diffeomorphism. Let $\psi\dvtx \R
\rightarrow(0,\infty)$
be its inverse; we have
%
%e12 #&#
%
\begin{eqnarray}\label{eqpsider}
\psi(z) & = & \theta^{-1}(z) = W\circ\exp(z+1),
\nonumber\\[-8pt]\\[-8pt]\nonumber
\psi^{(1)}(z) & = & \frac{1}{\theta^{(1)}\circ\psi(z)} = \frac{\psi(z)}{1+\psi
(z)}
\in(0,1),
\nonumber
\end{eqnarray}
where $W\dvtx (0,\infty)\rightarrow(0,\infty)$ is the Lambert $W$
function. In particular,
$\psi$ is strictly increasing, convex, and satisfies a linear growth
condition. So,
for any $\atil\in\Gtil$,
\[
\Emu\psi(\atil)^\lambda< K\bigl(1+\Emu\llvert \atil\rrvert
^\lambda\bigr) < \infty\quad\mbox{and}\quad \Emu\bigl\llvert \log\psi(
\atil)\bigr\rrvert ^\lambda= \Emu\bigl\llvert \atil-\psi(\atil)\bigr\rrvert
^\lambda < \infty.
\]
Let $P$ be the measure on $\cx$ with density $p=\psi(\atil)$; then
$P$ satisfies
(M1)--(M3), and $\tilde{\phi}(P)=\atil$, and this completes the proof.
\end{pf}

This construction defines an infinite-dimensional manifold of measures,
$(\Mtil,\Gtil,\tilde{\phi})$, with an atlas comprising the single
chart $\tilde{\phi}$. The
inverse map $\tilde{\phi}^{-1}\dvtx \Gtil\rightarrow\Mtil$ takes the form
%
%e13 #&#
%
\begin{equation}
\frac{\mathrm{d}\tilde{\phi}^{-1}(\atil)}{\mathrm{d}\mu}(x) = \psi\bigl(\atil(x)\bigr), \label{eqphitinv}
\end{equation}
where $\psi$ is as defined in (\ref{eqpsider}). (The definition
of $\psi$ used
here is slightly different from that in \cite{newt4}; in fact
$\psi_{\mathrm{here}}(z)=\psi_{\mathrm{there}}(z+1)$. The definition used here
has the advantage
that $\tilde{\phi}(\mu)=0$.) As in Section~\ref{seinfgeo}, we
consider a tangent vector
$U$ at $P\in\Mtil$ to be\vspace*{1pt} an equivalence class of differentiable
curves at $P$:  two
curves, $(\bfatil(t)\in\Gtil, t\in(-\varepsilon,\varepsilon))$ and
$(\bfbtil(t)\in\Gtil, t\in(-\varepsilon,\varepsilon))$, being
equivalent\vspace*{1pt} at $P$ if
$\bfatil(0)=\bfbtil(0)=\tilde{\phi}(P)$ and $\dot{\bfatil
}(0)=\dot{\bfbtil}(0)$. We denote
the tangent space at $P$ by $T_P\Mtil$, and the tangent bundle by
$T\Mtil:= \bigcup_{P\in\Mtil}(P,T_P\Mtil)$. The latter admits the
global chart
$\Phitil\dvtx T\Mtil\rightarrow\Gtil\times\Gtil$ where
%
%e14 #&#
%
\begin{equation}
\Phitil(P,U) = \bigl(\bfatil(0),\dot{\bfatil}(0)\bigr), \label{eqbunchartt}
\end{equation}
and $\bfatil$ is\vspace*{2pt} any differentiable curve in the equivalence class
$U$. If
$f\dvtx \Mtil\rightarrow Y$ is a map with range $Y$ (a Banach space) and
the map
$f\circ\tilde{\phi}^{-1}\dvtx \Gtil\rightarrow Y$ is (Fr\'{e}chet)
differentiable, then we write
\[
Uf:= \frac{\mathrm{d}}{\mathrm{d}t}\bigl(f\circ\tilde{\phi}^{-1}\bigr) \bigl(
\bfatil (t)\bigr)\biggm|_{t=0} = D\bigl(f\circ\tilde{\phi}^{-1}
\bigr)_\atil\util,
\]
where $(\atil,\util)=\Phitil(P,U)=(\tilde{\phi}(P),U\tilde{\phi})$.

%re3.1 #&#
%
\begin{remark} \label{reddiff}
The weaker notion of \textit{$d$-differentiability} is defined in \cite
{newt4}. In the
present context, the map $f\dvtx \Mtil\rightarrow Y$ is $d$-differentiable\vspace*{1pt}
if, for any
$P\in\Mtil$, there exists a continuous linear map
$d(f\circ\tilde{\phi}^{-1})_\atil\dvtx \Gtil\rightarrow Y$ such that
\[
\frac{\mathrm{d}}{\mathrm{d}t}\bigl(f\circ\tilde{\phi}^{-1}\bigr) \bigl(\bfatil(t)
\bigr)\biggm|_{t=0} = d\bigl(f\circ\tilde{\phi}^{-1}
\bigr)_\atil\util,
\]
for all differentiable curves $\bfatil$ in the equivalence class $U$.
(See Definition~3.1 in \cite{newt4}.) We then write $Uf=d(f\circ\tilde{\phi
}^{-1})_\atil\util$. Clearly,
if $f$ is Fr\'{e}chet differentiable then it is also
$d$-differentiable, and the
derivatives are identical. However, the converse is not always true, as
demonstrated
by Example 3.1 in \cite{newt4}.
\end{remark}

For any $\alpha\in[-1,1]$, let $\xi_\alpha\dvtx \R\rightarrow\R$ be
defined by
%
%e15 #&#
%
\begin{equation}
\xi_\alpha(z) = \cases{ \displaystyle \frac{2}{1-\alpha}
\psi(z)^{(1-\alpha)/2}, &\quad if $\alpha\in[-1,1)$,
\vspace*{3pt}\cr
\displaystyle \log\psi(z)
= z + 1 - \psi(z), &\quad if $\alpha=1$,} \label{eqxialpdef}
\end{equation}
where $\psi$ is as defined in (\ref{eqpsider}). Let $\xi_\alpha
^{(i)}$ be the $i$th
derivative of $\xi_\alpha$. An induction argument shows that all such
derivatives are
bounded, the first being
%
%e16 #&#
%
\begin{equation}
\xi_\alpha^{(1)}(z) = \frac{\psi(z)^{(1-\alpha)/2}}{1+\psi(z)} \in(0,1).
\label{eqxialpder}
\end{equation}
For any $\alpha\in[-1,1]$ and any $r\in[1,\lambda]$, let ${\Xi
_\alpha^r}\dvtx \Gtil
\rightarrow{L^r(\mu)}$ be
defined by
%
%e17 #&#
%
\begin{equation}
{\Xi_\alpha^r}(\atil) (x) = \xi_\alpha\bigl(\atil(x)
\bigr). \label{eqXialpdef}
\end{equation}
${\Xi_\alpha^r}$ is the \emph{superposition} (\emph{Nemytskij}) \emph{operator} associated
with the nonlinear
function $\xi_\alpha$, the domain $\Gtil$ and the range ${L^r(\mu
)}$. The
differentiability
properties of such operators are developed in an abstract setting in
Chapter~3 of
\cite{apza1}. In the present context, we are able to exploit the
explicit nature of
$\xi_\alpha$ to give a direct, self-contained proof of the following.

%le3.1 #&#
%
\begin{lem} \label{leXialpder}
%\begin{enumerate}
\textup{(i)}~${\Xi_\alpha^r}$ is of\vspace*{1pt} class $C^{\lceil\lambda/r\rceil
-1}$, where
$\lceil y\rceil:=\min\{i\in\bZ\dvt  y\le i\}$. For any $1\le i\le
\lceil\lambda/r\rceil-1$,
$D^i{\Xi_\alpha^r}\dvtx \Gtil\rightarrow L(\Gtil^i; {L^r(\mu)})$ is
given by
%
%e18 #&#
%
\begin{equation}
D^i{\Xi_{\alpha,\atil}^r}(\util_1,\ldots,
\util_i) (x) = {\xi _\alpha^{(i)}}\bigl(\atil(x)
\bigr)\util _1(x)\cdots\util_i(x). \label{eqXialpder}
\end{equation}

\textup{(ii)} ${\Xi_\alpha^r}$ satisfies global Lipschitz continuity
and linear
growth conditions,
and, for any $1\le i\le\lceil\lambda/r\rceil-1$,
%
%e19 #&#
%
\begin{equation}
\sup_{\atil\in\Gtil}\bigl\llVert D^i{\Xi_{\alpha,\atil}^r}
\bigr\rrVert < \infty. \label{eqXialpdb}
\end{equation}
%
%\end{enumerate}
%
\end{lem}

\begin{pf}
Let $l:=\lceil\lambda/r\rceil-1$, let $\atil,\util_1,\ldots,\util
_l\in\Gtil$, and let
$(\atil_n\neq\atil, n\in\N)$ be a sequence converging to $\atil$
in $\Gtil$. For
convenience of notation, let $\xi_\alpha^{(0)}:=\xi_\alpha$. If
$l\ge1$ then the mean value
theorem applied on an $x$-by-$x$ basis shows that, for any $0\le i\le l-1$,
\[
\bigl({\xi_\alpha^{(i)}}(\atil_n)-{
\xi_\alpha^{(i)}}(\atil)\bigr)\util _1\cdots
\util_i = \xi_\alpha^{(i+1)}(\atil)\util_1
\cdots\util_i(\atil_n-\atil) + R_n,
\]
where $R_n:= S_n\util_1\cdots\util_i(\atil_n-\atil)$, and, for some
$\beta=\beta(i,\atil_n(x),\atil(x))\in[0,1]$,
\[
S_n:= \xi_\alpha^{(i+1)}\bigl((1-\beta)\atil+\beta
\atil_n\bigr)-\xi _\alpha^{(i+1)}(\atil).
\]
Now $r(i+1)<\lambda$ and so, setting $s=\lambda/(\lambda-r(i+1))$,
we have
$1/s+(i+1)r/\lambda=1$ and H\"{o}lder's inequality shows that
\[
\Emu\llvert R_n\rrvert ^r \le\bigl(\Emu\llvert
S_n\rrvert ^{rs}\bigr)^{1/s}\bigl(\Emu\llvert
\util_1\rrvert ^\lambda\bigr)^{r/\lambda} \cdots\bigl(\Emu
\llvert \util_i\rrvert ^\lambda\bigr)^{r/\lambda}\bigl(\Emu
\llvert \atil_n-\atil \rrvert ^\lambda\bigr)^{r/\lambda}.
\]
So
\[
\llVert \atil_n-\atil\rrVert ^{-1}\sup
_{\llVert  \util_k\rrVert  =1}\llVert R_n\rrVert _{L^r(\mu)}\le\llVert
S_n\rrVert _{L^{rs}(\mu)}.
\]
Now $S_n\rightarrow0$ in probability, and is bounded by $2\sup_z\llvert  \xi
_\alpha^{(i+1)}(z)\rrvert  $,
and so it follows from the bounded convergence theorem that
$\llVert  S_n\rrVert  _{L^{rs}(\mu)}\rightarrow0$. An induction argument on $i$
thus establishes
that ${\Xi_\alpha^r}$ admits Fr\'{e}chet derivatives up to order $l$,
and that
these derivatives
take the form (\ref{eqXialpder}).

For any $0\le i\le l$, let
$T_n:= ({\xi_\alpha^{(i)}}(\atil_n)-{\xi_\alpha^{(i)}}(\atil
))\util_1\cdots\util_i$.
A similar argument
to that used above shows that
%
%e20 #&#
%
\begin{equation}
\sup_{\llVert  \util_k\rrVert  =1}\llVert T_n\rrVert _{L^r(\mu)} \le
\bigl\llVert {\xi_\alpha^{(i)}}(\atil_n)-{
\xi_\alpha^{(i)}}(\atil)\bigr\rrVert _{L^{rt}(\mu)},
\label{eqTbnd}
\end{equation}
where $t=\lambda/(\lambda-ri)$. If $i=0$ then the mean value theorem
and Jensen's
inequality show that
\[
\bigl\llVert {\xi_\alpha^{(i)}}(\atil_n)-{
\xi_\alpha^{(i)}}(\atil)\bigr\rrVert _{L^{rt}(\mu)} \le2\sup
_z\bigl\llvert \xi_\alpha^{(i+1)}(z)\bigr
\rrvert \llVert \atil_n-\atil\rrVert _\Gtil,
\]
which shows that ${\Xi_\alpha^r}$ satisfies global Lipschitz
continuity and
linear growth
conditions. If $i>0$ then $\llvert  {\xi_\alpha^{(i)}}(\atil_n)-{\xi_\alpha
^{(i)}}(\atil
)\rrvert  \rightarrow0$ in
probability, and is bounded by $2\sup_z\llvert  \xi_\alpha^{(i)}(z)\rrvert  $. In
either case, the
right-hand side of (\ref{eqTbnd}) converges to zero, and this shows
that ${\Xi_\alpha^r}$,
and any derivatives it has, are continuous. This completes the proof of
part (i).
The global boundedness of the derivatives~in~(\ref{eqXialpder})
follows from the
boundedness of those of $\xi_\alpha$, and this completes the proof of
part (ii).
\end{pf}

Lemma \ref{leXialpder} will be used in the following sections. It
also determines the
differentiability properties of Amari's $\alpha$-embedding maps,
$\Falp$, \cite{amar1}.
In the present context $\Falp\dvtx \Mtil\rightarrow L^2(\mu)$, and is
defined by
%
%e21 #&#
%
\begin{equation}
\Falp(P) = \Xi_\alpha^2\bigl(\tilde{\phi}(P)\bigr),
\label{eqFalpdef}
\end{equation}
where ${\Xi_\alpha^r}$ is as defined in (\ref{eqXialpdef}). The
choice of
$L^2(\mu)$ for
the range of $\Falp$ is consistent with the latter's role in the
definition of the
$\alpha$-divergences. (See (\ref{eql2rep}) and the expressions for
the derivatives of
$\cdal$ in Section~\ref{sestatdiv}.)

%co3.1 #&#
%
\begin{corol} \label{coFalpder}
For any $\alpha\in[-1,1]$, the map $\Falp$ is of class $C^{\lceil
\lambda
/2\rceil-1}$.
\end{corol}

\begin{pf}
This is an immediate consequence of Lemma \ref{leXialpder} with $r=2$.
\end{pf}

In the case $\lambda=2$, where $\Gtil$ is a Hilbert space, $\Falp$
is continuous but not
differentiable. It is, however, $d$-differentiable in the sense
described in
Remark \ref{reddiff}. (See Proposition 3.1 in \cite{newt4}.) More
generally, if
$\lambda=2n$ for some $n\in\N$, then $\Falp$ is of class $C^{n-1}$,
and its highest
Fr\'{e}chet derivative is $d$-differentiable. However, the
$d$-derivative may not be
continuous.

%s4 #&#
\section{The manifolds of probability measures} \label{semanifm}

Let $M$ be the subset of $\Mtil$ whose members are \emph{probability} measures.
These satisfy (M1)--(M3) and the additional hypothesis:
\begin{enumerate}
\item[(M4)] $\Emu p = 1$.
\end{enumerate}
Let $G$ ($=G_\lambda:=L_0^\lambda(\mu)$) be the Lebesgue space of
(equivalence classes
of) random variables $a\dvtx \bX\rightarrow\R$ for which $\Emu
\llvert  a\rrvert  ^\lambda<\infty$ and
$\Emu a=0$, and let $\phi\dvtx M\rightarrow G$ be defined as follows:
%
%e22 #&#
%
\begin{equation}
\phi(P) = \tilde{\phi}(P) - \Emu\tilde{\phi}(P) = p -1 + \log p - \Emu\log p.
\label{eqphidef}
\end{equation}

%pr4.1 #&#
%
\begin{prop} \label{prbijec}
\begin{longlist}[(iii)]
\item[(i)] $\phi$ is a bijection onto $G$.
\item[(ii)] $(M,G,\phi)$ is a $C^{\lceil\lambda\rceil-1}$-embedded
submanifold of
$(\Mtil,\Gtil,\tilde{\phi})$.
\item[(iii)] Let $\rho:=\tilde{\phi}\circ\phi^{-1}$ be the
inclusion map
$\imath\dvtx M\rightarrow\Mtil$ expressed in terms of the charts $\tilde
{\phi}$ and~$\phi$.
For any bounded set $B\subset G$,
%
%e23 #&#
%
\begin{equation}
\sup_{a\in B} \bigl(\bigl\llVert \rho(a)\bigr\rrVert +\llVert D
\rho_a\rrVert +\cdots+\bigl\llVert D^{\lceil
\lambda\rceil-1}
\rho_a\bigr\rrVert \bigr) < \infty. \label{eqrhobnd}
\end{equation}
\end{longlist}
\end{prop}

\begin{pf}
Let $l:=\lceil\lambda\rceil-1$. Let $\Psi\dvtx G\times\R\rightarrow
(0,\infty)$ be defined
by
\[
\Psi(a,z) = \Emu\psi(a+z) = \Emu\Xi_{-1}^1(a+z),
\]
where $\psi$ is as in (\ref{eqpsider}) and $\Xi_\alpha^r$ is as in
(\ref{eqXialpdef}). It follows from Lemma \ref{leXialpder}, with
$r=1$, that
$\Psi$ is of class $C^l$ and that, for any $u\in G$ and any $y\in\R$,
%
%e24 #&#
%
\begin{equation}
D\Psi_{a,z}(u,y) = \Emu{\psi^{(1)}}(a+z)u + \Emu{\psi
^{(1)}}(a+z)y. \label{eqPsider}
\end{equation}
For any $a\in G$, let $\Psi_a\dvtx \R\rightarrow(0,\infty)$ be defined by
$\Psi_a(z)=\Psi(a,z)$; then
\[
\Psi_a^{(1)}(z) = \Emu{\psi^{(1)}}(a+z) > 0.
\]
Since $\psi$ is convex,
\[
\sup_z\Psi_a \ge\sup_z
\psi\bigl(\Emu(a+z)\bigr) = \sup_z\psi(z) = +\infty;
\]
furthermore, the monotone convergence theorem shows that
\[
\lim_{z\downarrow-\infty} \Psi_a = \Emu\lim_{z\downarrow-\infty
}
\psi(a+z) = 0.
\]
Thus $\Psi_a\dvtx \R\rightarrow(0,\infty)$ is a bijection with strictly
positive derivative,
and the inverse function theorem shows that it is a $C^l$-isomorphism.
The implicit
function theorem shows that $Z\dvtx G\rightarrow\R$, defined by $Z(a) =
\Psi_a^{-1}(1)$,
is of class $C^l$. According to (\ref{eqPsider}), its first
derivative takes the form:
%
%e25 #&#
%
\begin{equation}
DZ_au = -\frac{\Emu{\psi^{(1)}}(a+Z(a))u}{\Emu{\psi
^{(1)}}(a+Z(a))}. \label{eqZder}
\end{equation}
Let $P$ be the probability measure on $\cx$ with density $p=\psi
(a+Z(a))$; then it
follows from (\ref{eqpsider}) and the mean value theorem that, for
any $x\in\bX$,
\[
\bigl\llvert p(x)-\psi\bigl(Z(a)\bigr)\bigr\rrvert \le\bigl\llvert a(x)\bigr
\rrvert\quad\mbox{and}\quad\bigl\llvert \log p(x)-\log\psi \bigl(Z(a)\bigr)\bigr
\rrvert \le\bigl\llvert a(x)\bigr\rrvert.
\]
So $P\in M$, and
\[
\phi(P) = \theta\circ\psi\bigl(a+Z(a)\bigr) - \Emu\log\psi\bigl(a+Z(a)\bigr) = a
+ Z(a) -\Emu\log\psi\bigl(a+Z(a)\bigr).
\]
Now $\phi(P)-a\in G$, and so
%
%e26 #&#
%
\begin{equation}
Z(a) = \Emu\log\psi\bigl(a+Z(a)\bigr) = - \cdpo(P\mid \mu), \label{eqZrep}
\end{equation}
and $\phi(P)=a$, which completes the proof of part (i).

According to (\ref{eqphidef}) and (\ref{eqZrep}), for any $a\in G$,
%
%e27 #&#
%
\begin{equation}
\rho(a) = a + \Emu\log\psi\bigl(a+Z(a)\bigr) = a + Z(a), \label{eqrhorep}
\end{equation}
and so $\rho$ is also of class $C^l$. $\rho$ is injective, as is its first
derivative; in fact, for any $\btil\in\rho(G)$ and any $\vtil\in
D\rho_aG$,
\[
\rho^{-1}(\btil) = \btil- \Emu\btil\quad\mbox{and}\quad D
\rho_a^{-1}\vtil= \vtil- \Emu\vtil,
\]
from which it also follows that $\rho$ and $D\rho_a$ are topological
embeddings.
Since $D\rho_a$ is also a linear map, it is a \emph{toplinear
isomorphism}, and its
image $D\rho_aG$ is a closed linear subspace of $\Gtil$. Suppose, in
the special
case that $\lambda=2$, that $\vtil\in\Gtil$ ($=\Gtil_2$) is
orthogonal to
$D\rho_aG$ ($=D\rho_aG_2$). It is a consequence of the representation
(\ref{eqZder}) that
\[
\Emu\vtil \biggl(u - \frac{\Emu{\psi^{(1)}}(\atil)u}{\Emu{\psi
^{(1)}}(\atil
)} \biggr) = \Emu\vtil D
\rho_au =0 \qquad\mbox{for all }u\in G_2,
\]
where $\atil=\rho(a)$. It then readily follows that
\[
\bigl\langle\vtil\Emu{\psi^{(1)}}(\atil)-{\psi^{(1)}}(\atil )
\Emu\vtil, u \bigr\rangle_{G_2} = 0\qquad\mbox{for all }u\in
G_2.
\]
So the orthogonal complement of $D\rho_aG_2$ in $\Gtil_2$ is the
one-dimensional
subspace,
%
%e28 #&#
%
\begin{equation}
E_a:= \bigl\{y{\psi^{(1)}}(\atil), y\in\R \bigr\}.
\label{eqEdef}
\end{equation}
Since ${\psi^{(1)}}$ is bounded, $E_a$ is also a one-dimensional
subspace of
$\Gtil_\lambda$
for any $\lambda\in[2,\infty)$. Now $D\rho_aG_\lambda=\Gtil
_\lambda\cap D\rho_aG_2$,
and so $D\rho_aG_\lambda\oplus E_a = \Gtil_\lambda$ and
$D\rho_aG_\lambda\cap E_a = \{0\}$.
We have thus shown that $D\rho_a$ \emph{splits} $\Gtil$ into the
complementary closed
subspaces $D\rho_aG$ and $E_a$. It thus follows from Proposition 2.3
of Chapter II in
\cite{lang1} that $\rho$ is a {\em$C^l$-immersion}. Since $\rho$ is
a topological
embedding it is also a $C^l$-embedding, and this completes the proof of
part (ii).

It follows from Jensen's inequality and (\ref{eqpsider}) that, for
any $a\in G$,
%
%e29 #&#
%
\begin{eqnarray}
-\log\Emu{\psi^{(1)}}(\atil) & \le& -\Emu\log{\psi^{(1)}}(
\atil)
\nonumber
\\
& = & -\Emu\log p + \Emu\log(1+p) \label{eqlogrho}
\\
& \le& \cdpo(P\mid \mu) + \log2,
\nonumber
\end{eqnarray}
where $\atil=\rho(a)$ and $P=\phi^{-1}(a)$. Now
%
%e30 #&#
%
\begin{eqnarray}
\cdmo(P\mid \mu)+\cdpo(P\mid \mu) & = & \Emu(p-1) \bigl(\log p+\cdpo(P\mid \mu)
\bigr)
\nonumber
\\
& \le& \Emu \bigl(p-1+\log p +\cdpo(P\mid \mu) \bigr)^2/2
\label{eqrelbnd}
\\
& \le& \llVert a\rrVert _G^2/2,
\nonumber
\end{eqnarray}
and so, since they are both non-negative, $\cdmo(\phi^{-1}\mid \mu)$ and
$\cdpo(\phi^{-1}\mid \mu)$ are bounded on bounded sets. Together with
(\ref{eqlogrho}),
this proves that
%
%e31 #&#
%
\begin{equation}
\inf_{a\in B} \Emu{\psi^{(1)}}\bigl(\rho(a)\bigr) > 0.
\label{eqdenbnd}
\end{equation}

The boundedness of the derivatives of $\rho$ on bounded sets follows from
(\ref{eqdenbnd}), the boundedness of the derivatives of $\psi$, and
an induction
argument. The boundedness of $\rho$ on bounded sets follows from (\ref
{eqZrep}),
(\ref{eqrhorep}) and (\ref{eqrelbnd}), and this completes the proof
of part (iii).
\end{pf}

The tangent space at base point $P\in M$, $T_PM$, can be defined in the
same way as
was $T_P\Mtil$. The tangent bundle, $TM:= \bigcup_{P\in M}(P,T_PM)$,
admits the global
chart $\Phi\dvtx TM\rightarrow G\times G$, where
%
%e32 #&#
%
\begin{equation}
\Phi(P,U) = \bigl(\bfa(0),\dot{\bfa}(0)\bigr) = \bigl(\phi(P),U\phi\bigr),
\label
{eqbunchart}
\end{equation}
and $\bfa$ is any differentiable curve in the equivalence class $U$.
For any $P\in M$,
the tangent space $T_PM$ is a subspace of $T_P\Mtil$ of co-dimension
1; in fact
%
%e33 #&#
%
\begin{equation}
T_P\Mtil= T_PM \oplus\{yU_0, y\in\R\},
\label{eqtansplit}
\end{equation}
where $U_0$ is the equivalence class of differentiable curves on $\Mtil
$ containing
the curve $(\bfatil(t):=\atil+t{\psi^{(1)}}(\atil), t\in
(-\varepsilon
,\varepsilon))$, and
$\atil=\tilde{\phi}(P)$. (See (\ref{eqEdef}).)

%co4.1 #&#
%
\begin{corol} \label{cofaider}
The map $\Falp\circ\imath\dvtx M\rightarrow L^2(\mu)$, where $\Falp$ is
as defined in
(\ref{eqFalpdef}), is of class $C^{\lceil\lambda/2\rceil-1}$.
\end{corol}

\begin{pf}
This follows from Corollary \ref{coFalpder} and Proposition \ref
{prbijec}(ii).
\end{pf}

%s5 #&#
\section{The \texorpdfstring{$\alpha$}{$alpha$}-divergences} \label{sestatdiv}

We begin by investigating the regularity of the $\alpha$-divergences
on $\Mtil$. The
usual extension of the $\alpha$-divergences of (\ref{eqalphdiv}) to
sets of finite
measures such as $\Mtil$ is as follows \cite{amar1}:
%
%e34 #&#
%
\begin{equation}
\cdal(P\mid Q):= \cases{ Q(\bX) - P(\bX) + \Emu p\log(p/q),\qquad\mbox{if $\alpha=-1$},
\vspace*{3pt}\cr
\displaystyle\frac{2}{1+\alpha}P(\bX) + \frac{2}{1-\alpha}Q(\bX) -
\frac{4}{1-\alpha^2}\Emu p^{(1-\alpha)/2}q^{(1+\alpha)/2},
\vspace*{3pt}\cr
\hspace*{151pt}\mbox{if $\alpha \in(-1,1)$},
\vspace*{3pt}\cr
P(\bX) - Q(\bX) + \Emu q\log(q/p),\hspace*{20,5pt}\mbox{if $\alpha=1$.}}
\label{eqalphdivt}
\end{equation}
These can be represented in terms of the maps $\Falp$ of (\ref{eqFalpdef});
for example,
%
%e35 #&#
%
\begin{equation}
\frac{4}{1-\alpha^2}\Emu p^{(1-\alpha)/2}q^{(1+\alpha)/2} = \bigl\langle\Falp(P),
\Fmalp(Q) \bigr\rangle_{L^2(\mu)}. \label{eql2rep}
\end{equation}
So, for any $\alpha\in[-1,1]$ and any $P,Q\in\Mtil$, $\cdal(P\mid
Q) <
\infty$,
and we refer to elements of $\Mtil$ as ``finite-entropy'' measures. We could
investigate the smoothness properties of $\cdal$ starting\vspace*{1pt} from those
of $\Falp$.
However, this approach would show only that the divergences are of class
$C^{\lceil\lambda/2\rceil-1}$; a stronger result can be obtained by
a more direct
approach. The following lemma, which is similar in nature to Lemma \ref
{leXialpder},
prepares the ground. For any $\alpha\in[-1,1]$, let
$\Upalp\dvtx \Gtil\times\Gtil\rightarrow L^1(\mu)$ be the following
superposition operator:
%
%e36 #&#
%
\begin{equation}
\Upalp(\atil,\btil) (x) = \xi_\alpha\bigl(\atil(x)\bigr)\ximalp\bigl(
\btil(x)\bigr), \label{eqUpalpdef}
\end{equation}
where $\xi_\alpha$ is as in (\ref{eqxialpdef}).

%le5.1 #&#
%
\begin{lem} \label{leUpalpder}
For any $0\le i,j\le\lfloor\lambda\rfloor-1$ with $i+j\le\lceil
\lambda\rceil-1$, the
map $\Upalp$ is of class $C^{i,j}$. Its partial derivatives,
$\Upalp^{(i,j)}:= D_1^iD_2^j\Upalp\dvtx \Gtil\times\Gtil\rightarrow
L(\Gtil
^{i+j};L^1(\mu))$, are
given by
%
%e37 #&#
%
\begin{eqnarray}\label{eqUpalpij}
&& \Upalp^{(i,j)}(\atil,\btil) (\util_1,\ldots,
\util_i;\vtil_1\ldots,\vtil_j) (x) \qquad
\nonumber\\[-8pt]\\[-8pt]\nonumber
&&\quad = {\xi_\alpha^{(i)}}\bigl(\atil(x)\bigr){
\ximalp^{(j)}}\bigl(\btil (x)\bigr)\util _1(x)\cdots
\util_i(x) \vtil_1(x)\cdots\vtil_j(x).
\end{eqnarray}
\end{lem}

\begin{pf}
Let $0\le i\le\lfloor\lambda\rfloor-2$, $0\le j\le\lfloor\lambda
\rfloor-1$ and
$i+j\le\lceil\lambda\rceil-2$. Let
$\atil,\btil$, $\util_1,\ldots,  \util_i,\vtil_1,\ldots,\vtil_j\in
\Gtil$ and let
$(\atil_n\neq\atil, n\in\N)$ be a sequence converging to $\atil$
in $\Gtil$. Applying
the mean value theorem on an $x$-by-$x$ basis, we obtain
\[
\bigl({\xi_\alpha^{(i)}}(\atil_n)-{
\xi_\alpha^{(i)}}(\atil)\bigr){\ximalp ^{(j)}}(\btil)
\util_1\cdots \vtil_j = \xi_\alpha^{(i+1)}(
\atil){\ximalp^{(j)}}(\btil)\util_1\cdots \vtil
_j(\atil_n-\atil) + R_n,
\]
where $R_n:= S_n\util_1\cdots\vtil_j(\atil_n-\atil)$ and, for some
$\beta=\beta(i,j,\atil_n(x),\atil(x),\btil(x))\in[0,1]$,
\[
S_n:= \bigl(\xi_\alpha^{(i+1)}\bigl((1-\beta)\atil+
\beta\atil_n\bigr) -\xi_\alpha^{(i+1)}(\atil) \bigr){
\ximalp^{(j)}}(\btil).
\]
Now $i+j+1<\lambda$ and so, setting $s=\lambda/(\lambda-i-j-1)$, we have
$1/s+(i+j+1)/\lambda=1$, and H\"{o}lder's inequality shows that
\[
\Emu\llvert R_n\rrvert \le \bigl(\Emu\llvert S_n\rrvert
^s\bigr)^{1/s}\bigl(\Emu\llvert \util_1\rrvert
^\lambda\bigr)^{1/\lambda} \cdots\bigl(\Emu\llvert \vtil_j
\rrvert ^\lambda\bigr)^{1/\lambda}\bigl(\Emu\llvert \atil_n-
\atil \rrvert ^\lambda\bigr)^{1/\lambda},
\]
so that
\[
\llVert \atil_n-\atil\rrVert ^{-1}\sup
_{\llVert  \util_k\rrVert  =\llVert  \vtil_k\rrVert  =1}\llVert R_n\rrVert _{L^1(\mu)} \le\llVert
S_n\rrVert _{L^s(\mu)}.
\]
Now $S_n\rightarrow0$ in probability, and is dominated by
$f:=2\sup_z\llvert  \xi_\alpha^{(i+1)}(z)\rrvert  \llvert  {\ximalp^{(j)}}(\btil)\rrvert  $. If $j=0$
then $f\in\Gtil$ and
$s\le\lambda$, whereas if $j\ge1$ then $f\in L^\infty(\mu)$ and
$s<\infty$. In
either case the dominated convergence theorem shows that
$\llVert  S_n\rrVert  _{L^s(\mu)}\rightarrow0$, so that $\Upalp^{(i,j)}$ is
differentiable in
its first argument, with derivative $\Upalp^{(i+1,j)}$.
Similarly, if $0\le i\le\lfloor\lambda\rfloor-1$, $0\le j\le
\lfloor\lambda\rfloor-2$,
and $i+j\le\lceil\lambda\rceil-2$ then $\Upalp^{(i,j)}$ is
differentiable in
its second argument, with derivative $\Upalp^{(i,j+1)}$. An induction argument
on $i$ and $j$ thus establishes (\ref{eqUpalpij}).

It remains to show that, for any $0\le i,j\le\lfloor\lambda\rfloor
-1$ with
$i+j=\lceil\lambda\rceil-1$, $\Upalp^{(i,j)}$ is continuous. Now
\[
\bigl(\Upalp^{(i,j)}(\atil_n,\btil_n)-
\Upalp^{(i,j)}(\atil,\btil )\bigr) (\util_1,\ldots,
\vtil_j) = (T_{1,n}+T_{2,n}+T_{3,n})
\util_1\cdots\vtil_j,
\]
where
\begin{eqnarray*}
T_{1,n} &:= & \bigl({\xi_\alpha^{(i)}}(
\atil_n)-{\xi_\alpha ^{(i)}}(\atil)\bigr){
\ximalp^{(j)}}(\btil), \qquad T_{2,n}:= {\xi_\alpha^{(i)}}(
\atil) \bigl({\ximalp^{(j)}}(\btil _n)-{\ximalp^{(j)}}(
\btil)\bigr),
\\
T_{3,n} &:= & \bigl({\xi_\alpha^{(i)}}(
\atil_n)-{\xi_\alpha ^{(i)}}(\atil)\bigr) \bigl({
\ximalp^{(j)}}(\btil _n)-{\ximalp^{(j)}}(\btil)
\bigr),
\end{eqnarray*}
and similar arguments to those used above show that
\[
\bigl\llVert \Upalp^{(i,j)}(\atil_n,\btil_n)-
\Upalp^{(i,j)}(\atil,\btil)\bigr\rrVert \le\llVert T_{1,n}\rrVert
_{L^t(\mu)} + \llVert T_{2,n}\rrVert _{L^t(\mu)} + \llVert
T_{3,n}\rrVert _{L^t(\mu)},
\]
where $t=\lambda/(\lambda-i-j)$. We will thus have established the
continuity of
$\Upalp^{(i,j)}$ if we can show that
%
%e38 #&#
%
\begin{equation}
\llVert T_{k,n}\rrVert _{L^t(\mu)} \rightarrow0\qquad\mbox{for }k=1,2,3. \label{eqTcon}
\end{equation}
If $i=j=0$, then $t=1$ and (\ref{eqTcon}) follows from the
Cauchy--Schwarz inequality
and the mean value theorem; for example,
\begin{eqnarray*}
\llVert T_{1,n}\rrVert _{L^1(\mu)}^2 & \le& \bigl
\llVert \xi_\alpha(\atil_n)-\xi_\alpha(\atil)\bigr
\rrVert _{L^2(\mu)}\bigl\llVert \ximalp(\btil)\bigr\rrVert _{L^2(\mu)}
\\
& \le& \sup_z\bigl\llvert \xi_\alpha^{(1)}(z)
\bigr\rrvert \llVert \atil_n-\atil\rrVert _{L^2(\mu)} \sup
_z\bigl\llvert \ximalp^{(1)}(z)\bigr\rrvert \llVert
\btil\rrVert _{L^2(\mu)} \rightarrow0.
\end{eqnarray*}
If $i,j>0$, then $t<\infty$, and $T_{k,n}\rightarrow0$ in probability
and is bounded
for all $k$; so (\ref{eqTcon}) follows from the bounded convergence
theorem. If
$i=0$ and $j>0$, then $t\le\lambda$ and
\[
\llVert T_{k,n}\rrVert _{L^t(\mu)} \le2\sup_z
\bigl\llvert {\ximalp^{(j)}}(z)\bigr\rrvert \bigl\llVert
\xi_\alpha(\atil_n)-\xi_\alpha (\atil )\bigr\rrVert
_{L^t(\mu)} \rightarrow0\qquad\mbox{for }k=1,3;
\]
furthermore $T_{2,n}\rightarrow0$ in probability and is dominated by
$2\sup_z\llvert  {\ximalp^{(j)}}(z)\rrvert  \llvert  \xi_\alpha(\atil)\rrvert  \in\Gtil$, and so the
dominated convergence theorem
establishes (\ref{eqTcon}). The case $i>0$ and $j=0$ can be treated
in the same way,
and this completes the proof.
\end{pf}

%co5.1 #&#
%
\begin{corol} \label{codivder}
For any $\alpha\in[-1,1]$, and any $0\le i,j\le\lfloor\lambda
\rfloor-1$ with
$i+j\le\lceil\lambda\rceil-1$, the \mbox{$\alpha$-}divergence
$\cdal\dvtx \Mtil\times\Mtil\rightarrow[0,\infty)$ is of class $C^{i,j}$.
\end{corol}

\begin{pf}
It follows from (\ref{eqalphdiv}), (\ref{eqxialpder}), (\ref
{eqXialpdef}) and
(\ref{eqUpalpdef}) that
\[
\cdal(P\mid Q):= \cases{ \displaystyle \Emu \bigl(\Xi_{-1}^1(
\btil)-\Xi_{-1}^1(\atil)+\Upsilon _{-1}(\atil,
\atil) -\Upsilon_{-1}(\atil,\btil) \bigr), &\quad if $\alpha=-1$,
\vspace*{3pt}\cr
\displaystyle\Emu \biggl(\frac{2}{1+\alpha}\Xi_{-1}^1(
\atil)+\frac{2}{1-\alpha
}\Xi_{-1}^1(\btil)- \Upalp(\atil,
\btil) \biggr), &\quad if $\alpha\in(-1,1)$,
\vspace*{3pt}\cr
\displaystyle\Emu \bigl(
\Xi_{-1}^1(\atil)-\Xi_{-1}^1(\btil)+
\Upsilon_{1}(\btil,\btil) -\Upsilon_{1}(\atil,\btil)
\bigr), &\quad if $\alpha=1$,}
\]
where $\atil=\tilde{\phi}(P)$ and $\btil=\tilde{\phi}(Q)$. The
corollary thus follows from Lemma
\ref{leXialpder} (with $r=1$) and Lemma \ref{leUpalpder}.
\end{pf}

Straightforward calculations show that, for any $\atil,\btil,\util,\vtil\in\Gtil$,
%
%e39 #&#
%
\begin{eqnarray}\label{eqdivder1}
D_1\cdal\bigl(\tilde{\phi}^{-1}\mid \tilde{
\phi}^{-1}\bigr)_{\atil,\btil
}\util & = & \Emu\bigl(
\Upalp^{(1,0)}(\atil,\atil)-\Upalp^{(1,0)}(\atil,\btil )\bigr)\util,
\nonumber\\[-8pt]\\[-8pt]\nonumber
D_2\cdal\bigl(\tilde{\phi}^{-1}\mid \tilde{
\phi}^{-1}\bigr)_{\atil,\btil
}\vtil & = & \Emu\bigl(
\Upalp^{(0,1)}(\btil,\btil)-\Upalp^{(0,1)}(\atil,\btil )\bigr)\vtil.
\end{eqnarray}
If $\lambda>2$, these admit the following representations
%
%e40 #&#
%
\begin{eqnarray}\label{eqbdivder1}
U\cdal(\fndot\mid Q) & = & \bigl\langle\Fmalp(P)-\Fmalp(Q), U\Falp\bigr
\rangle_{L^2(\mu)},
\nonumber\\[-8pt]\\[-8pt]\nonumber
V\cdal(P \mid\fndot) & = & \bigl\langle\Falp(Q)-\Falp(P), V\Fmalp\bigr
\rangle_{L^2(\mu)},
\nonumber
\end{eqnarray}
and $\cdal$ admits the following mixed second derivative
%
%e41 #&#
%
\begin{equation}
D_1D_2\cdal\bigl(\tilde{\phi}^{-1}\mid
\tilde{\phi}^{-1}\bigr)_{\atil,\btil
}(\util,\vtil) = - \Emu
\Upalp^{(1,1)}(\atil,\btil) (\util,\vtil) = - \langle U\Falp, V\Fmalp
\rangle_{L^2(\mu)}, \label{eqdivder2}
\end{equation}
where $(P,U)=\Phitil^{-1}(\atil,\util)$ and $(Q,V)=\Phitil
^{-1}(\btil,\vtil)$.
Setting $\btil=\atil$, we obtain the following definition of the
(extended) Fisher
metric on $T_P\Mtil$:  for any $U,V\in T_P\Mtil$,
%
%e42 #&#
%
\begin{equation}
\langle U, V \rangle_P:= -UV\cdal= \langle U\Falp, V\Fmalp
\rangle_{L^2(\mu)}. \label{eqfisher}
\end{equation}

%re5.1 #&#
%
\begin{remark} \label{redivddiff}
The representations in (\ref{eqbdivder1}) and (\ref{eqdivder2}),
and the definition
in (\ref{eqfisher}), are also valid for the case $\lambda=2$ if the
weaker notion of
$d$-differentiability is used in the definitions of $U\Falp$, $V\Fmalp
$ and
$UV\cdal$. (See \cite{newt4}.)
\end{remark}

It follows from (\ref{eqxialpder}) and (\ref{eqdivder2}) that
$\langle V, U \rangle_P = \langle U, V \rangle_P$, and that, for any
$s\in\R$,
$\langle sU, V \rangle_P=\langle U, sV \rangle_P = s\langle U, V
\rangle_P$.
Furthermore,
%
%e43 #&#
%
\begin{equation}
\langle U, U \rangle_P = \Emu\frac{p}{(1+p)^2}\util^2
\le\Emu\util^2 \le\llVert \util\rrVert _\Gtil^2,
\label{eqfishdom}
\end{equation}
where $\util=U\tilde{\phi}$; in particular $\langle U,U\rangle_P=0$
if and only
if $\util=0$. Thus, $(T_P\Mtil,\langle\fndot\rangle_P)$ is an
inner product
space. As shown in (\ref{eqfishdom}), the Fisher norm is dominated by
the natural
Banach norm on $T_P\Mtil$. However, it is not equivalent to that norm,
even in the
case $\lambda=2$. (See \cite{newt4}.) In the general,
infinite-dimensional case
$(T_P\Mtil,\langle\fndot\rangle_P)$ is not a Hilbert space; the
Fisher metric
is a pseudo-Riemannian metric but not a Riemannian metric.

If $\lambda>3$, $\cdal$ admits the following mixed third derivative
%
%e44 #&#
%
\begin{equation}
D_1^2D_2\cdal\bigl(\tilde{
\phi}^{-1}\mid\tilde{\phi}^{-1}\bigr)_{\atil,\btil
}(\util,
\vtil;\wtil) = -\Emu\Upalp^{(2,1)}(\atil,\btil) (\util,\vtil;\wtil)
\label
{eqdivder3}.
\end{equation}
Setting $\btil=\atil$ and carrying out some straightforward
calculations, we obtain
%
%e45 #&#
%
\begin{equation}
D_1^2D_2\cdal\bigl(\tilde{
\phi}^{-1}\mid\tilde{\phi}^{-1}\bigr)_{\atil,\atil
}(\util,
\vtil;\wtil) = -\Emu\frac{p}{(1+p)^2}\Gamalptil(\atil,\util,\vtil)\wtil,
\label{eqbdivder3}
\end{equation}
where $\Gamalptil\dvtx \Gtil\times\Gtil\times\Gtil\rightarrow
L^{\lambda/2}(\mu)$ is defined by
%
%e46 #&#
%
\begin{equation}
\Gamalptil(\atil,\util,\vtil) (x) = \frac{1-\alpha}{2}\frac{\util(x)\vtil(x)}{(1+p(x))^2} -
\frac{1+\alpha}{2}p(x)\frac{\util(x)\vtil(x)}{(1+p(x))^2}. \label{eqGamalptdef}
\end{equation}
If $\atil=\tilde{\phi}(P)$, and $\util$ and $\vtil$ are such that
$\Gamalptil(\atil,\util,\vtil)\in\Gtil$, then there exist tangent vectors
$Y,W\in T_P\Mtil$ such that $\Gamalptil(\atil,\util,\vtil)=Y\tilde
{\phi}$ and
$\wtil=W\tilde{\phi}$. In this case
%
%e47 #&#
%
\begin{equation}
D_1^2D_2\cdal\bigl(\tilde{
\phi}^{-1}\mid\tilde{\phi}^{-1}\bigr)_{\atil,\atil
}(\util,
\vtil;\wtil) = -\langle Y, W \rangle_P. \label{eqfishrep}
\end{equation}

For any $l\in\N_0$ and any $s\in[1,\infty]$, let $\cvtil_s^l$ be
the set of vector
fields $\bfV\dvtx \Mtil\rightarrow T\Mtil$ for which
$\bfvt(P) (:=\bfV(P)\tilde{\phi})\in L^{s\lambda}(\mu)$ for all
$P\in\Mtil$, and $\bfvt\dvtx \Mtil\rightarrow L^{s\lambda}(\mu)$ is of
class $C^l$. For
any $\bfU\in\cvtil_s^0$, we can use (\ref{eqfishrep}) and the
Eguchi relations
\cite{eguc1} to define an ``$\alpha$-derivative''
$\nablatil_\bfU^\alpha\dvtx \cvtil_{s/(s-1)}^1\rightarrow\cvtil_1^0$,
as follows
%
%e48 #&#
%
\begin{equation}
\nablatil_\bfU^\alpha\bfV:= \Phitil^{-1} \bigl(
\tilde{\phi},\bfU\bfvt+ \Gamalptil(\tilde {\phi},\bfut,\bfvt) \bigr).
\label{eqalpder}
\end{equation}
However, this does not define an operator, $\nablatil^\alpha$, with domain
$\cvtil_1^0\times\cvtil_1^1$, and so it does not define a full covariant
derivative on $T\Mtil$. With the exception of the $+1$ connection on
the exponential
Orlicz manifold, this appears to be an insuperable problem in infinite
dimensions.
In order for the divergences to be sufficiently smooth, the tangent
space must be
given a stronger topology than that generated by the Fisher metric, and
so it is
incomplete with respect to the latter. This creates difficulties with
the projection
methods at the heart of the definition of $\alpha$-covariant
derivatives. In the
special case that $s=\infty$, $\nablatil_\bfU^\alpha\bfV$ is well
defined for all
$C^1$ vector fields $\bfV$, and thus provides a limited notion of
$\alpha$-parallel
transport on the tangent bundle. (See \cite{gras1} for a similar
result on the
exponential Orlicz manifold.)

A straightforward calculation shows that, for any $\alpha\in[-1,1]$,
$\bfU
\in\cvtil_s^0$ and
$\bfV,\bfW\in\cvtil_{s/(s-1)}^1$,
%
%e49 #&#
%
\begin{equation}
\bfU\langle\bfV, \bfW\rangle_P = \bigl\langle\nablatil_\bfU^\alpha
\bfV, \bfW\bigr\rangle_P + \bigl\langle\bfV, \nablatil_\bfU^{-\alpha}
\bfW\bigr\rangle_P, \label{eqcodual}
\end{equation}
reflecting the duality (\ref{eqexpcodu}) of the finite-dimensional case.

%s5.1 #&#
\subsection{The Fisher metric and \texorpdfstring{$\alpha$}{$alpha$}-derivatives on \texorpdfstring{$(M,G,\phi)$}{$(M,G,phi)$}} \label{sealphm}

In the above, we used the $\alpha$-divergences and Eguchi relations to
define the
extended Fisher metric and $\alpha$-derivatives on the manifold $\Mtil
$. Clearly, we
could follow the same approach with the submanifold $M$. (It follows
from Proposition
\ref{prbijec}(ii) and Corollary \ref{codivder} that the $\alpha
$-divergences have
the same smoothness properties on $M$ as they have on $\Mtil$.) For
any $P\in M$,
the Fisher metric on $T_PM$, thus obtained, is a restriction of the
extended Fisher
metric on $T_P\Mtil$, as defined above. (See (\ref{eqtansplit}).) On
the other hand,
the definition of the $\alpha$-derivative involves second derivatives
of $\cdal$ in one variable, and so the correspondence between $\Mtil$
and $M$ is not so transparent.

For some $s\in[1,\infty]$, let $\bfU\in\cvtil_s^0$ and $\bfV\in
\cvtil_{s/(s-1)}^1$ be
vector fields on $\Mtil$, whose restrictions to $M$ are vector fields
of $M$; then,
for any $P\in M$, $\Phitil(\bfU(P))=(\rho(a),D\rho_a{\mathbf u}(P))$ and
$\Phitil(\bfV(P))=(\rho(a),D\rho_a{\mathbf v}(P))$, where $(a,{\mathbf
u}(P))=\Phi(\bfU(P))$ and
$(a,{\mathbf v}(P))=\Phi(\bfV(P))$. So, according to (\ref{eqalpder}),
%
%e50 #&#
%
\begin{eqnarray}
\nablatil_\bfU^\alpha\bfV & = & \Phitil^{-1} \bigl(
\tilde{\phi},\bfU(D\rho{\mathbf v}) + \Gamalptil(\tilde{\phi},D\rho{\mathbf u},D\rho{\mathbf
v}) \bigr)
\nonumber
\\
& = & \Phitil^{-1} \bigl(\tilde{\phi},D\rho\bfU{\mathbf v}+
D^2\rho ({\mathbf u},{\mathbf v}) + \Gamalptil(\tilde{\phi},D\rho{\mathbf u},D
\rho{\mathbf v}) \bigr)
\\
& = & \Phitil^{-1} \biggl(\tilde{\phi},D\rho\bfU{\mathbf v}+
\frac
{1-\alpha}{2}\gamma - \frac{1+\alpha}{2}\eta \biggr),
\nonumber
\end{eqnarray}
where $\gamma,\eta\dvtx \Mtil\rightarrow\Gtil$ are defined by
%
%e51 #&#
%
\begin{eqnarray}\label{eqgameta}
\gamma(P) (x) & = & \frac{D\rho_a{\mathbf u}(P)(x)D\rho_a{\mathbf v}(P)(x)}{(1+p(x))^2} + D^2\rho_a
\bigl({\mathbf u}(P),{\mathbf v}(P)\bigr) (x),
\nonumber\\[-8pt]\\[-8pt]
\eta(P) (x) & = & p(x)\frac{D\rho_a{\mathbf u}(P)(x)D\rho_a{\mathbf v}(P)(x)}{(1+p(x))^2} - D^2
\rho_a\bigl({\mathbf u}(P),{\mathbf v}(P)\bigr) (x),
\nonumber
\end{eqnarray}
and $a=\phi(P)$. It follows from (\ref{eqZder}) and (\ref
{eqrhorep}) that, for any
$u\in G$,
\[
D\rho_au = u - \frac{\Emu{\psi^{(1)}}(\rho(a))u}{\Emu{\psi
^{(1)}}(\rho(a))} = u - \frac{\Emu D\Xi_{-1,\rho(a)}^1u}{\Emu D\Xi_{-1,\rho(a)}^11},
\]
and so, according to the quotient and chain rules of differentiation,
and Lemma
\ref{leXialpder},
\[
D^2\rho_a(u,v) = -\frac{\Emu\psi^{(2)}(\rho(a))D\rho_au D\rho_av}{\Emu{\psi^{(1)}}
(\rho(a))} = -
\frac{1}{\Emu{\psi^{(1)}}(\rho(a))}\Emu{\psi^{(1)}}\bigl(\rho (a)\bigr)\frac{D\rho
_au D\rho_av}{(1+p)^2}.
\]
From these derivatives and (\ref{eqgameta}), we conclude that
\begin{eqnarray*}
\gamma(P) & = & D\rho_a \biggl(\frac{\bfut(P)\bfvt(P)}{(1+p)^2} - \Emu
\frac{\bfut(P)\bfvt(P)}{(1+p)^2} \biggr),
\\
\eta(P) & = & D\rho_a \biggl(p\frac{\bfut(P)\bfvt(P)}{(1+p)^2} - p\bigl\langle
\bfU(P), \bfV(P) \bigr\rangle_P \biggr)
\\
& &{} + \bigl\langle\bfU(P), \bfV(P) \bigr\rangle_P \biggl(p-
\frac{\Emu{\psi^{(1)}}(\atil)p}{\Emu{\psi^{(1)}}(\atil
)} \biggr)
\\
& & {}+ \frac{1}{\Emu{\psi^{(1)}}(\atil)}\Emu{\psi ^{(1)}}(\atil) \biggl(p
\frac{\bfut(P)\bfvt(P)}{(1+p)^2} +\frac{\bfut(P)\bfvt(P)}{(1+p)^2} \biggr)
\\
& = & D\rho_a \biggl(p\frac{\bfut(P)\bfvt(P)}{(1+p)^2} - p\bigl\langle\bfU(P),
\bfV(P) \bigr\rangle_P \biggr)
\\
& & {}+ \bigl\langle\bfU(P), \bfV(P) \bigr\rangle_P \biggl(p-
\frac{1}{\Emu{\psi^{(1)}}(\atil)}+1 \biggr) + \frac
{1}{\Emu
{\psi^{(1)}}(\atil)}\Emu p \frac{\bfut(P)\bfvt(P)}{(1+p)^2}
\\
& = & D\rho_a \biggl(p\frac{\bfut(P)\bfvt(P)}{(1+p)^2} - p\bigl\langle\bfU(P),
\bfV(P) \bigr\rangle_P \biggr) + (1+p)\bigl\langle\bfU(P), \bfV(P)
\bigr\rangle_P,
\end{eqnarray*}
where $\atil=\tilde{\phi}(P)$, $\bfut(P)=D\rho_a{\mathbf u}(P)$,
$\bfvt(P)=D\rho_a{\mathbf v}(P)$, and we
have used the fact that ${\psi^{(1)}}\psi=\psi-{\psi^{(1)}}$ in the
second step.
We have thus shown
that
\[
\nablatil_\bfU^\alpha\bfV(P) = \Phitil^{-1} \biggl(
\atil,D\rho_a \bigl(\bfU{\mathbf v}(P) + \Gamalp \bigl(a,{\mathbf u}(P),{\mathbf
v}(P)\bigr) \bigr) - \frac{1+\alpha}{2}(1+p)\bigl\langle\bfU(P), \bfV(P) \bigr
\rangle_P \biggr),
\]
where $\Gamalp\dvtx G\times G\times G\rightarrow L_0^{\lambda/2}(\mu)$ is
defined by
%
%e52 #&#
%
\begin{eqnarray}\label{eqGamalpdef}
\Gamalp(a,u,v) (x) & = & \frac{1-\alpha}{2} \biggl(\frac{D\rho_au(x)D\rho_av(x)}{(1+p(x))^2} - \Emu
\frac{D\rho_au(x)D\rho_av(x)}{(1+p(x))^2} \biggr)
\nonumber\\[-8pt]\\[-8pt]\nonumber
& & {}- \frac{1+\alpha}{2}p(x) \biggl(\frac{D\rho_au(x)D\rho
_av(x)}{(1+p(x))^2} - \langle
U, V \rangle_P \biggr).
\nonumber
\end{eqnarray}
For any $W\in T_PM$
%
%e53 #&#
%
\begin{eqnarray}
\bigl\langle\nablatil_\bfU^\alpha\bfV(P), W \bigr
\rangle_P & = & \Emu\frac{p}{(1+p)^2}D\rho_a \bigl(
\bfU{\mathbf v}(P) + \Gamalp\bigl(a,{\mathbf u}(P),{\mathbf v}(P)\bigr) \bigr)D
\rho_aw
\nonumber
\\
&&{}- \frac{1+\alpha}{2}\bigl\langle\bfU(P), \bfV(P) \bigr\rangle
_P\Emu\frac{p}{(1+p)^2} (1+p)D\rho_aw \label{eqprojcd}
\\
& = & \bigl\langle\nabla_\bfU^\alpha\bfV(P), W \bigr
\rangle_P,
\nonumber
\end{eqnarray}
where $\nabla_\bfU^\alpha\bfV\dvtx M\rightarrow TM$ is the vector field
on $M$ defined by
%
%e54 #&#
%
\begin{equation}
\nabla_\bfU^\alpha\bfV = \Phi^{-1} \bigl(\phi,\bfU{
\mathbf v}+ \Gamalp(\phi,{\mathbf u},{\mathbf v}) \bigr). \label{eqnablalpdef}
\end{equation}
As (\ref{eqprojcd}) shows, $\nabla_\bfU^\alpha\bfV(P)$ is the
projection of
$\nablatil_\bfU^\alpha\bfV(P)$ onto $T_PM$, in the sense of the
Fisher metric.
The map $\nabla_\bfU^\alpha\dvtx \cv_{s/(s-1)}^1\rightarrow\cv_1^0$,
thus defined, is the
$\alpha$-derivative on $M$, which could also found by direct
calculation in the same way
as was $\nablatil_\bfU^\alpha$.

%s6 #&#
\section{Finite-dimensional submanifolds} \label{sefindim}

For some $d\in\N$ and $n\in\N\cup\{\infty\}$, let $(N,B,\theta)$
be a $d$-dimensional
$C^n$-embedded submanifold of $\Mtil$. By this, we mean that $N\subset
\Mtil$, $B$ is
a non-empty open subset of $\R^d$, $\theta\dvtx N\rightarrow B$ is a
bijection, and the
inclusion map $\imathtil\dvtx N\rightarrow\Mtil$ is both a topological
embedding and a
$C^n$-immersion. (See, e.g., \cite{lang1}.) As in Section~\ref
{seinfgeo}, the
tangent space at base point $P\in N$, $T_PN$, is spanned by the vectors
$(\partial_i, i=1,\ldots,d)$, where $\partial_i$ is the equivalence
class of
differentiable curves on $N$ containing the curve
$(\bfy_i(t):=\theta(P)+t\bfe_i, t\in(-\varepsilon,\varepsilon))$. The
matrix form of the
(extended) Fisher metric is
%
%e55 #&#
%
\begin{equation}
g(P)_{i,j}:= \langle\partial_i,\partial_j
\rangle_P = \Emu\frac{p}{(1+p)^2}\wtil_i
\wtil_j, \label{eqfishmat5}
\end{equation}
where $\wtil_i=\partial_i\tilde{\phi}$.

Since $(T_PN, \langle\fndot\rangle_P)$ is a \emph{finite-dimensional}
inner-product space it is also a Euclidean space, and the Fisher metric
is a
Reimannian metric on $N$. If $\lambda>3$ and $n\ge2$, the full theory of
$\alpha$-covariant derivatives and their associated geometries can
thus be developed
on $N$. According to the Eguchi relations, the Christoffel symbols for the
$\alpha$-covariant derivative on $(N,\theta)$ are
%
%e56 #&#
%
\begin{equation}
\Gamma_\alpha^N(P)_{i,j}^k:=
-g(P)^{k,l}\partial_i\partial_j
\partial_l\cdal = g(P)^{k,l}\Emu\frac{p}{(1+p)^2}\Gamalptil
\bigl(\tilde{\phi}(P),\wtil _i,\wtil_j\bigr)
\wtil_l, \label{eqchrissym5}
\end{equation}
where $g(P)^{k,l}$ is the $(k,l)$ element of the inverse of the matrix $g(P)$,
$\partial_i$ and $\partial_j$ act on the first argument of $\cdal$,
$\partial_l$ acts on the second argument of $\cdal$, and $\Gamalptil
$ is as
defined in~(\ref{eqGamalptdef}).

If $N$ is a \emph{statistical} manifold (it is also a subset of $M$)
then the
inclusion map, $\imath\dvtx N\rightarrow M$, takes the form $\imath=\pi
\circ\imathtil$,
where $\pi=\phi^{-1}\circ\rhobar\circ\tilde{\phi}$ and $\rhobar
\dvtx \Gtil\rightarrow G$ is
defined by $\rhobar(\atil)=\atil-\Emu\atil$. Clearly $\pi$ is of
class $C^\infty$,
and so $\imath$ is of class $C^n$. Furthermore, $\partial_i\imathtil
\in T_PM$ for
all $i$, and the restriction of the pushforward $\pi_*$ to $T_PM$ is
the identity map
of $T_PM$, and so the derivative of $\imath$ is injective. Since
$\imathtil$ is a
topological embedding and the map $\rho$ of Proposition \ref
{prbijec} is continuous,
$\imath$ is also a topological embedding. It thus follows that $N$ is
also a
$C^n$-embedded submanifold of $M$.

We finish with two examples of finite-dimensional submanifolds that
illustrate the
foregoing developments.

%ex6.1 #&#
%
\begin{example} \label{exphitsub}
Let $\eta_1,\ldots, \eta_d$ be linearly independent elements of
$\Gtil$, let
$\gamma\dvtx \R^d\rightarrow\Gtil$ be defined by $\gamma(y) = y^i \eta
_i$, and let
$N:=\tilde{\phi}^{-1}\circ\gamma(\R^d)$. Since\vspace*{1pt} the $\eta_i$ are linearly
independent $\gamma$ is an injection, and $(N,\R^d,\theta)$, with
$\theta:=\gamma^{-1}\circ\phi$, is a $d$-dimensional manifold. It
is trivially a
$C^\infty$-embedded submanifold of $\Mtil$.
\end{example}

%ex6.2 #&#
%
\begin{example} \label{exexpsub}
Let $(N,B,\theta)$ be the $d$-dimensional exponential statistical
manifold defined
in Section~\ref{seinfgeo}, where the underlying space $(\bX,\cx,\mu
)$ coincides with
that of Sections~\ref{semanifmt}--\ref{sestatdiv}, and suppose that
the $\eta_i$
and $B$ are such that $\theta^{-1}(B)\subseteq M$. It is shown in
Theorem 5.1 of
\cite{newt4} that $N$, thus defined, is a $C^\infty$-embedded
submanifold of $M$.
(Strictly speaking, Theorem 5.1 in \cite{newt4} addresses only the
case $\lambda=2$;
however, the same proof carries over to the more general setting where
$\lambda\in[2,\infty)$.)
\end{example}

%s7 #&#
\section{Concluding remarks}

Because of their role in the definition of the Kullback--Leibler
divergence, it is
natural to regard the density, $p$, and its log as belonging to dual
function spaces. The choice of the exponential Orlicz space for $\log
p$ (and,
implicitly, its dual for $p$) yields the manifold of \cite{pise1},
comprising all
probability measures in an absolute continuity equivalence class. The
choice in
\cite{newt4} of the Hilbert space $L_0^2(\mu)$ for both $p$ and $\log
p$ leads to a
significantly simpler construction, but at a cost to inclusiveness.
This is also true
of the Banach space approach taken here. However, this is unimportant
in many
applications (and may even be beneficial). In problems of Bayesian
estimation, for
example, we do not need manifolds to contain more than the posterior
distributions
associated with the various observations, and some finite-dimensional
structures on
which approximations can be based.

The choice of reference measure $\mu$ is important. The use of a \emph{finite} measure
is natural in the context of (M2) and (M3), and since the elements and
topologies of
$\Mtil$ and $M$ are not affected by its total mass, it is also natural
to assume that
$\mu$ is a \emph{probability} measure. If $\bX=\R^n$ then \textup{(M1)} is
satisfied by all
measures that are mutually absolutely continuous with respect to
Lebesgue measure if,
for example, $\mu$ is a non-singular multi-variate Gaussian measure.
It may seem that
one could construct larger manifolds by piecing together coordinate patches
$(\Mtil_i,\Gtil_i,\tilde{\phi}_i,\mu_i)$, defined as in
Section~\ref{semanifmt} but with
differing patch-centric measures, $\mu_i$. However, this is not
possible since the
requirement that $\mathrm{d}P/\mathrm{d}\mu_i\in L^\lambda(\mu_i)$ for each $i$ is
incompatible with the
regularity of the transition maps $\tilde{\phi}_i\circ\tilde{\phi
}_j^{-1}$ in all but trivial
cases (such as that in which $\mathrm{d}\mu_i/\mathrm{d}\mu_j\in L^\infty(\mu_j)$ for
all $i,j$). The
requirement that $\mathrm{d}P/\mathrm{d}\mu\in L^\lambda(\mu)$ is stronger than needed
for pure information
geometry (even in the Hilbert case, $\lambda=2$). However, it is
useful in its own
right. For example, in the context of Bayesian estimation, it yields
bounds such as
(\ref{eql2bnd}).

The role played by the exponential function in an exponential family,
such as that
of \cite{pise1}, is played here by the function $\psi$. In this
sense, $\Mtil$ and
$M$ are extreme examples of \emph{general deformed families}, as defined
in Chapter~10
of \cite{naud1}. (They are extreme in the sense that $\psi$ satisfies
a linear growth
condition.) General deformed families of probability measures are also
developed and
generalised in \cite{vica1}, where they are referred to as \emph{$\varphi$-families}.
The function $\varphi$ is used there in the definition of the
(Musielak--Orlicz) model
spaces, and gives rise to dual divergence functions distinct from
$\cdal$. Here, our
aim is somewhat different from those of \cite{naud1} and \cite
{vica1}. We provide a
simple framework for the classical information geometry in infinite
dimensions; this
requires a stronger topology on the model space than that associated
with the
$\varphi$-function $\psi$. (The Musielak--Orlicz spaces associated
with the
$\varphi$-function $\psi$ have topologies that are too weak, even for
the definition
of the Fisher metric.)

\section*{Acknowledgement}

The author would like to thank an anonymous referee for suggesting the
offset $-1$ in
the definition of $\tilde{\phi}$ in (\ref{eqphitdef}), which
introduces a number of
advantages.

% imsref loaded by linak, 2014-11-17 09:15:55
% imsref loaded by linak, 2014-11-18 15:19:16
% imsref loaded by linak, 2014-11-18 15:20:52

\printhistory
\end{document}